\documentclass[12pt]{amsart}
\usepackage{amsmath,
amsfonts,
amsthm,
amssymb,
mathrsfs}
\usepackage{amsmath}
\usepackage{amsfonts}
\usepackage{amssymb}

\newcommand{\comment}[1]{}

\newcommand{\Hall}{\mathrm{Hull}}

\newcommand{\Seq}[1]{\langle #1 \rangle}

\newcommand{\ZFC}{\mathrm{ZFC}}

\newcommand{\Lsc}{\textsc{L}}

\newcommand{\supp}{\mathrm{supp}}
\newcommand{\rest}{\upharpoonright}

\newcommand{\pphi}{\varphi}

\newcommand{\otp}{\mathrm{otp}}

\theoremstyle{plain}
\newtheorem{theorem}{Theorem}[]
\newtheorem{lem}[theorem]{Lemma}
\newtheorem{prop}[theorem]{Proposition}

\newtheorem{question}[theorem]{Question}

\theoremstyle{definition}
\newtheorem{defn}[theorem]{Definition}

\begin{document}

\title{Aronszajn free  Kurepa trees}
\author[H. Lamei Ramandi, S. Todorcevic]{Hossein Lamei Ramandi and Stevo Todorcevic}
\address{Institute f\"{u}r Mathematische Logik und 
Grundlagenforschung  \\
Westf\"{a}lische Wilhelms-Universit\"{a}t 
M\"{u}nster,  Germany}
\email{hlamaira@exchange.wwu.de}

\address{Department of Mathematics \\ University of Toronto,
Toronto \\ Canada}
\email{{\tt stevo@math.toronto.edu}}

\address{ Institut de Math\'ematiques de Jussieu, Paris, France}
\email{{\tt stevo.todorcevic@imj-prg.fr}}

\address{ Matemati\v{c}ki Institut, SANU, Belgrade, Serbia}
\email{{\tt stevo.todorcevic@sanu.ac.rs}}

\subjclass[]{}

\begin{abstract} 
We consider a transitive relation on $\mathcal{P}(\omega_1)$ and show 
if there is a maximal element with respect to this relation then there is a Kurepa tree with no Aronszajn subtree.
We also show that if there is a maximal element in $\mathcal{P}(\omega_1)$, then there are Kurepa trees which are not 
club isomorphic.
These maximal subsets of $\omega_1$ exist 
in many known models that are obtained from 
the constructible universe $\Lsc$ without large cardinal assumptions. 
For instance, whenever 
$\alpha_0 \in \omega_1$ and $X \subset \omega_1$ are such that 
$\omega_1^{\textsc{L}[X \cap \alpha_0]} = \omega_1, \omega_2^{\Lsc[X]} = \omega_2$ 
and $\textsc{V}$ is a semiproper forcing 
extension of $\textsc{L}[X]$ then 
$X$ is maximal in $\textsc{V}$. 
\end{abstract}

\maketitle
\section{introduction}
Let us call an $\omega_1$-tree \emph{Aronszajn free} if it has no Aronszajn subtree.
We are interested in the question that under what circumstances 
all Kurepa trees contain Aronszajn subtrees?
Jensen showed that there is an Aronszajn free Kurepa 
tree in the construtible universe $\Lsc$ (see \cite{Aronszajn_free_Kurepa_of_Jensen}).
Todorcevic showed that there is a countably closed forcing notion which 
adds an Aronszajn free Kurepa tree (\cite{trees_subtrees: Todorcevic}).
In \cite{Komjath_Aronszajn_Kurepa}, 
assuming there are two inaccessible cardinals, 
Komjath proves the consistency 
of the statement that every Kurepa tree has an Aronszajn subtree. 
The large cardinal assumption for the consistency of  this statement 
is reduced to existence of one inaccessible cardinal
in \cite{Aronszajn_Kurepa}. 
In other words, if there is an 
inaccessible cardinal then it is consistent that 
every Kurepa tree contains an Aronszajn subtree.
It is not known if we need an inaccessible cardinal 
in order to obtain the consistency of every Kurepa tree has an 
Aronszajn subtree.
Motivated by this background, we show the following theorem.
\begin{theorem}\label{main}
Assume there are 
 $\alpha_0 \in \omega_1$ and  $ X \subset \omega_1$ such that: 
 \begin{enumerate}
 \item \label{LC}$\omega_1^{\Lsc[X \cap \alpha_0]} = \omega_1$,
  $\omega_2^{\Lsc[X]} = \omega_2$, and 
   \item \label{semiproper} $\textsc{V}$ is a semiproper forcing 
   extension of $\Lsc[X]$. 
 \end{enumerate}
Then there is an Aronszajn free Kurepa tree.
\end{theorem}

Note  that if there is no inaccessible cardinal in $\Lsc$ then there are 
$\alpha_0 \in \omega_1$ and $X \subset \omega_1$
for which  (\ref{LC}) of Theorem \ref{main} holds.
In particular, there are no semiproper forcing notion 
which adds Aronszajn subtrees to Kurepa trees of $\Lsc$
and which keeps them Kurepa.
For more motivation and history 
behind finding Aronszajn free Kurepa trees in  
various models of $\ZFC$ we refer the reader to 
\cite{Aronszajn_free_Kurepa_of_Jensen}, \cite{Sikorski} and \cite{trees:Todorcevic}.

In order to prove Theorem \ref{main}, 
we consider a transitive relation  on $\mathcal{P}(\omega_1)$ -- the power set of $\omega_1$. 
The motivation behind this relation can be explained as follows.
Whenever 
there are maximal elements with respect to this relation
and there are no inaccessible cardinals in $\Lsc$ 
then there is an Aronszajn free Kurepa tree. 
Moreover, if $X$ satisfies the conditions (\ref{LC}) and (\ref{semiproper}) of Theorem \ref{main} then $X$ is maximal.
This is more clarified in Definition \ref{order} and Propositions \ref{Tree} and \ref{max_in_semi}.

Maximal subsets of $\omega_1$ can possibly be useful in various 
constructions regarding problems around the 
second order properties of objects of size $\aleph_1$. 
For instance, a similar construction to Proposition \ref{Tree} 
can be used in order to obtain two Kurepa trees which are not club isomorphic. For more clarity, we include this construction in Proposition \ref{Similar} which together with Proposition \ref{max_in_semi}
implies  the following theorem. 

\begin{theorem} \label{club_iso_theorem}
Assume there are 
 $\alpha_0 \in \omega_1$ and  $ X \subset \omega_1$ such that: 
 \begin{enumerate}
 \item $\omega_1^{\Lsc[X \cap \alpha_0]} = \omega_1$,
  $\omega_2^{\Lsc[X]} = \omega_2$, and 
   \item $\textsc{V}$ is a semiproper forcing 
   extension of $\Lsc[X]$. 
 \end{enumerate}
Then there are Kurepa trees that are not club isomorphic.
\end{theorem}

For Proposition \ref{Tree}, we were inspired by Jensen's construction which appears in \cite{Aronszajn_free_Kurepa_of_Jensen}.
The advantage of our construction can be explained as follows.
In the process of diagonalization of our construction 
we do not need to refer to ``potential Aronszajn subtrees''.
This is a key to our construction because 
the process of diagonalization can happen in some $\Lsc$-like model and we do not 
need to deal with $\mathcal{P}^{\textsc{V}}(\omega_1)$.
This is not the case in Jensen's Construction, where 
in the process of diagonalization one has to consider all possible 
Aronszajn subtrees and carry out the diagonalization against 
Aronszajn subtrees as well as $\omega_1$-many branches.

\section{proof of theorems}

In order to avoid ambiguity, we fix some notations. 
If $\alpha$ is an ordinal and $t:A \longrightarrow \alpha$ is a function, support of $t$, which is denoted by $\supp (t)$, is the set of all elements $a \in A$ with $t(a) \neq 0$. 
If $T$ is a tree and $\alpha$ is an ordinal, 
$T(\alpha)$ is the $\alpha$'th level of $T$.
If $T$ is a tree and $A$ is a set of ordinals, 
$T \rest A$ is the tree consisting of all elements of $T$ whose height 
is in $A$. The order on $T \rest A$ is the inherited order from $T$.

Assume $A$ is a set of ordinals, $\kappa$ is a regular cardinal in $\Lsc[A]$  and 
$x \in \Lsc_\kappa[A]$. $\Hall(x;\Lsc_\kappa[A])$ is the smallest countable $M \prec \Lsc_\kappa [A]$ with $x \in M$.
This notation makes sense because $\Lsc[X]$ has a canonical 
well ordering and Skolem functions are well defined.
The following Lemma, in a slightly different form, can be found in 
\cite{Aronszajn_free_Kurepa_of_Jensen}. We include its proof for more 
clarity.
\begin{lem}[\cite{Aronszajn_free_Kurepa_of_Jensen}] \label{magic}
Assume $X \subset \omega_1$, 
$\kappa > \omega_1$ is a regular cardinal, 
$x \in \Lsc_\kappa [X]$, 
$\alpha = \omega_1 \cap \Hall (x; \Lsc_\kappa [X])$, 
$\omega_1 = \omega_1^{\Lsc[X]}$ and 
$\gamma = \otp (\kappa \cap \Hall (x;\Lsc_\kappa[X]))$.
Then $\alpha$ is countable in $\Lsc_{\gamma+2}[X \cap \alpha]$.
\end{lem}
\begin{proof}
Let $\pphi : \Hall(x; \Lsc_{\kappa}[X]) \longrightarrow \Lsc_{\gamma}[X \cap \alpha]$ be the transitive collapse map.
Assume without loss of generality that $\alpha$ is uncountable in 
$\Lsc_{\gamma+1}[X \cap \alpha]$.
So in particular, $\alpha$ is the first uncountable 
ordinal in $\Lsc_{\gamma+1}[X \cap \alpha]$.
For every $n \in \omega$, let $Y_n$ be the smallest $Y \prec_{\Sigma_n} \Lsc_{\gamma}[X \cap \alpha]$
such that $\pphi(x) \in Y$ and $Y \cap \alpha$ is transitive.
Note that $Y_n$ is definable in $\Lsc_{\gamma}[X \cap \alpha]$
and in particular, $Y_n \in \Lsc_{\gamma+1}[X \cap \alpha]$.
 
Let $\alpha_n = Y_n \cap \alpha$.
Note that $\alpha _n \in \alpha$ since $Y_n$ is an element and countable in  $\Lsc_{\gamma+1}[X \cap \alpha]$. 
Since $\Seq{Y_n : n \in \omega}$ is definable in 
$\Lsc_{\gamma+1}[X \cap \alpha]$, so is 
$\Seq{\alpha_n : n \in \omega}$.
Since $\bigcup_{n \in \omega} Y_n = \Lsc_{\gamma}[X \cap \alpha]$, 
$\Seq{\alpha_n : n \in \omega}$ is a cofinal sequence of $\alpha$
which is in $\Lsc_{\gamma+2}[X \cap \alpha]$, as desired.
\end{proof}

\begin{defn}\label{order}
Assume $X \subset \omega_1$, $\kappa \geq \omega_2$ is a regular cardinal and $\beta \in \kappa$. Define $h_\kappa(X,\beta) = \omega_1 \cap \Hall(\{X, \beta \} ; \Lsc_\kappa[X])$.
Let $X <_{\kappa} Y$ if for all $\beta \in \kappa$, 
$h_\kappa(X, \beta) < h_\kappa (Y, \beta)$. 
A set $X \subset \omega_1$ is said to be $\kappa$-maximal
if there is no $Y\subset \omega_1$ with 
$X<_\kappa Y$.
\end{defn}

Note that 
for every regular $\kappa \geq \omega_2$
the empty set is $\kappa$-maximal  in 
$\mathcal{P}(\omega_1)^{\Lsc}$. 
Similarly
$X$ is $\kappa$-maximal in $\mathcal{P}(\omega_1)^{\Lsc[X]}$, for every regular  $\kappa \geq \omega_2$. 
\begin{prop}\label{Tree}
Assume $\kappa \geq 2^{\omega_1}$ is a regular cardinal, $\alpha_0 \in \omega_1, X \subset \omega_1$,
$X$ is $\kappa$-maximal, 
$\omega_1^{\textsc{L}[X \cap \alpha_0]} = \omega_1,$ and 
$ \omega_2^{\Lsc[X]} = \omega_2$ .
Then there is a Kurepa tree with no Aronszajn subtree.
\end{prop}
\begin{proof}
Since $2^{< \omega_2}$ does not add new subsets of $\omega_1$, 
without loss of generality we can assume that $2^{\omega_1} = \omega_2$.
We define $T(\alpha) \subset 2^{\alpha}$ and $B_\alpha$ for $\alpha \in \omega_1$ inductively as follows.
For $\alpha \in \alpha _0 +1$ let $T(\alpha) = \{t \in 2^\alpha: \supp(t) \textrm{ is finite.}\} $ and $B_\alpha$ be the empty set.
Let $T(\alpha +1)= \{ t \in 2^{\alpha +1}: t \rest \alpha \in T(\alpha) \}$
and $B_{\alpha+1 } = \emptyset$
when $T(\alpha)$ and $ B_\alpha$ are given. 

Assume $\alpha \in \omega_1 \setminus \alpha_0$ is limit and 
$\Seq{T(\xi), B_\xi: \xi \in \alpha}$
are given such that 
$\bigcup_{\xi \in \alpha}T(\xi)$ is a downward closed subtree of $2^{\alpha}$.
Let $f(\alpha )$ be the smallest limit $\xi \in \omega_1$ such that $\alpha $ is countable in $\Lsc_\xi[X \cap \alpha]$. 
Let $B_\alpha $ be the collection of all $((U',U), \beta)$ in $(H_\kappa)^2 \times \kappa$ such that: 
\begin{itemize}
\item $U' \subset U \subset 2^{< \omega_1}$ are downward closed $\omega_1$-trees,
\item $U'$ is Aronszajn, 
\item $U \rest \alpha = \bigcup_{\xi \in \alpha} T(\xi)$,
\item $\omega_1 \cap \Hall(\{X,U, U', \beta \};\Lsc_\kappa[X, U, U']) = \alpha$ and 
\item $\omega_1 \cap \Hall(\{X,U, \beta \};\Lsc_\kappa[X, U]) = \alpha$.
\end{itemize}
For  a cofinal branch $b \subset \bigcup_{\xi \in \alpha}T(\xi)$, 
we allow $\bigcup b \in T(\alpha)$ if 
\begin{enumerate}
\item[$(e_1)$] $\bigcup b \in \Lsc_{f(\alpha)}[X \cap \alpha , \bigcup_{\xi \in \alpha} T(\xi)]$ and 
\item[$(e_2)$] there is no $((U', U), \beta) \in B_\alpha$ such that 
$\bigcup b$ is cofinal in $U'$.
\end{enumerate}
It is easy to see that 
$T$ is a downward closed subtree of $2^{< \omega_1}$.
Moreover, 
$\{ t \in 2^{< \omega_1} : \supp(t) \textrm{ is finite.} \}$ is a downward closed subtree of $T$. 
In particular, by $(e_1)$, $T$ is an uncountable $\omega_1$-tree.

In order to see
$T$ has no Aronszajn subtree, assume for a contradiction that 
$A \subset T$ is a downward closed Aronszajn subtree.
Since $X$ is $\kappa$-maximal, there is $\beta \in \kappa$ and $\delta \in \omega_1$
such that: 
\begin{itemize}
\item$\delta = \omega_1 \cap \Hall(\{X, T, \beta \}; \Lsc_\kappa [X, T])$ and 
\item $ \delta = \omega_1 \cap \Hall(\{X,A. T, \beta \}; \Lsc_\kappa [X,A, T])$.
\end{itemize}
Then $((A,T), \beta) \in B_\delta$ and in particular $A \cap T(\delta) = \emptyset$ by $(e_2)$. But this is a contradiction.

Now we show that $T$ is Kurepa. 
Assume for a contradiction that $\Seq{b_\xi : \xi \in \omega_1}$
is the $\Lsc[X,T]$-minimal  enumeration of all cofinal branches of $T$.
Define $\Seq{M_\alpha: \alpha \in \omega_1}$ inductively as follows:
\begin{itemize}
\item $M_0 = \Hall(\{X,T\};\Lsc_{\omega_2}[X,T])$, 
\item $M_{\alpha+1} = \Hall(\{ M_\alpha \}; \Lsc_{\omega_2}[X,T])$, and 
\item $M_\alpha = \bigcup_{\xi \in \alpha}M_\xi$ if $\alpha \in \omega_1$ is limit.
\end{itemize} 
Let $\delta_\xi = M_\xi \cap \omega_1$, $\gamma_\xi = \otp(M_\xi \cap \omega_2)$, 
and 
$C = \{ \delta_\xi : \xi \in \omega_1 \}$.
Note that $C \subset \omega_1$ is a club and $M_\xi$ is 
isomorphic to $\Lsc_{\gamma_\xi}[X \cap \delta_\xi , T \rest \delta_\xi]$. 
Moreover $C \cap \delta_\xi \in \Lsc_{\gamma_\xi + 3} [X \cap \delta_\xi , T \rest \delta_\xi]$,
by the same argument as in the proof of $\diamondsuit^+$ in $\Lsc$. 

In order to reach a contradiction, 
we define inductively an increasing sequence $\Seq{t_\alpha : \alpha \in C}$ in $T$ such that $\bigcup_{\alpha \in C}t_\alpha$ is different from all branches in  $\Seq{b_\xi : \xi \in \omega_1 }$. 
Let $t_{\delta_0}$ be the $\Lsc[X,T]$-minimal element in $T(\delta_0 +1)$ which is not in $\bigcup_{\xi \in \delta_0} b_\xi$.
If $t_{\delta_\xi}$ is given let $t_{\delta_{\xi +1 }}$ be the 
$\Lsc[X,T]$-minimal element in $T(\delta_{\xi +1} +1)$ which is above 
$t_{\delta_\xi}$ and which is not in $\bigcup_{\eta \in \delta_{\xi +1}} b_\eta.$
If $\xi \in \omega_1$ is limit, let $t_{\delta_\xi} =  \bigcup_{\eta \in C \cap \delta_\xi}t_\eta$.

We show by induction that $t_\eta$ is defined for every 
$\eta \in C$. 
The successor step is obvious. 
Assume $\xi$ is limit, $\delta = \delta_\xi$ and $\Seq{t_\eta: \eta \in C \cap \delta}$ is given.
Let $M = M_\xi$ and 
$\gamma = \gamma_\xi$.
Note that $\Seq{b_\xi : \xi \in \omega_1}$ is definable from $X,T$
and $C \cap \delta \in \Lsc_{\gamma+3}[X \cap \delta, T \rest \delta]$.
Therefore, $\bigcup_{\eta \in C \cap \delta} t_\eta$ is in 
$\Lsc_{f(\alpha)}[X \cap \delta , T \rest \delta]$, which shows 
$(e_1)$ holds for $\Seq{t_\eta: \eta \in C \cap \delta}$.

We need to show that $(e_2)$  holds for $\Seq{t_\eta: \eta \in C \cap \delta}$. 
Assume $((U',U),\beta) \in B_\delta$ and we will show that 
$\{t_\eta: \eta \in C \cap \delta \}$ is not a subset of $U'$.
Let $N = \Hall(\{X, U', U, \beta \};\Lsc_{\kappa}[X, U', U])$. 
Obviously, $\delta= N \cap \omega_1= \omega_1 \cap 
\Hall(\{X,U,\beta  \}; \Lsc_\kappa[X,U])$.
Let $\gamma' = \omega_1 \cap \Hall(\{X,U,\beta  \}; \Lsc_\kappa[X,U])$
and $\gamma'' = \otp(N \cap \kappa)$.
It is also obvious that $\gamma' \leq \gamma''$.

We show that
$\gamma < \gamma'$.
Recall that $U \rest \delta= T\rest \delta$.
Lemma \ref{magic} implies that $\delta$ is countable in 
$\Lsc_{\gamma' +2}[X \cap \delta , U \rest \delta]
= \Lsc_{\gamma' +2}[X \cap \delta , T \rest \delta]$.
On the other hand, $\delta$ is uncountable in  $\Lsc_{\gamma}[X \cap \delta , T \rest \delta]$, since this model is isomorphic to $M$.
But then $\gamma \leq \gamma'$, because $\gamma , \gamma'$ are 
both limit ordinals.
In order to see $\gamma' \neq \gamma$, observe that 
\begin{itemize}
\item[$(a)$] 
$\Hall(x, \Lsc_{\gamma}[X \cap \delta , T \rest \delta])$ is a proper subset of $\Lsc_{\gamma}[X \cap \delta , T \rest \delta]$, 
for all $x \in \Lsc_{\gamma}[X \cap \delta , T \rest \delta]$,  but 
\item[$(b)$] $\Lsc_{\gamma'}[X \cap \delta , T \rest \delta] = \Hall(x;\Lsc_{\gamma'}[X \cap \delta , T \rest \delta])$ for $x = \psi (\{X,U,\beta \})$ where $\psi:$
$\Hall(\{X,U,\beta  \}; \Lsc_\kappa[X,U]) \longrightarrow \Lsc_{\gamma'}[X \cap \delta , T \rest \delta]$ 
is the transitive collapse isomorphism.
\end{itemize}
Therefore, $\Lsc_{\gamma'}[X \cap \delta , T \rest \delta] \neq \Lsc_{\gamma}[X \cap \delta , T \rest \delta]$ and $\gamma \neq \gamma'$.

Let $t = \bigcup \{ t_\eta : \eta \in C \cap \delta \}$.
Recall that $t \in \Lsc_{\gamma +3}[X \cap \delta , T \rest \delta]$.
Then $t \in \Lsc_{\gamma'}[X \cap \delta , U \rest \delta]$
since $\gamma < \gamma'$ and $ T \rest \delta = U \rest \delta$.
Consequently, 
$t \in \Lsc_{\gamma''}[X \cap \delta , U \rest \delta, U' \rest \delta]$
since $\gamma'' \geq \gamma'$.
Let $\pphi$ denote the inverse of the transitive collapse map from   $\Lsc_{\gamma''}[X \cap \delta , U \rest \delta, U' \rest \delta]$ to $N$. 
Then $\pphi(t) \in N$ is a cofinal branch in $U$. 
Therefore $\pphi(t)$ is not a subset of $U'$.
Since $\pphi(U' \rest \delta) = U'$, $t$ is not a subset of $U'\rest \delta$, as desired. 
\end{proof}

\begin{prop}\label{Similar}
Assume $\kappa \geq 2^{\omega_1}$ is a regular cardinal, 
$\alpha_0 \in \omega_1$, 
$X \subset \omega_1$, 
$X$ is $\kappa$-maximal, 
$\omega_1^{\Lsc[X \cap \alpha_0]} = \omega_1$, and $\omega_2 = \omega_2^{\Lsc[X]}$.
Then there are Kurepa trees which are not  club isomorphic.
\end{prop}
\begin{proof}
We define two $\omega_1$-trees $S$ and $T$,  which are downward closed
subtrees of $2^{\omega_1}$, by simultaneous induction on the levels of these trees.
For $\alpha \in \alpha_0 +1$, let $T(\alpha) = S(\alpha) = \{ t \in 2^\alpha: 
\supp(t)$ is finite $\}$.
Assume $\alpha \in \lim(\omega_1) \setminus (\alpha_0 +1 ), T \rest \alpha :=
\bigcup_{\xi \in \alpha} T(\xi),$ and $ S \rest \alpha : = \bigcup_{\xi \in \alpha} S(\xi)$ are given. Define 
$f(\alpha)$ to be the smallest  
$  \xi \in \lim(\omega_1) \setminus (\alpha+1)$ such that 
every ordinal is countable in $\Lsc_{\xi}[X \cap \alpha ]$.
If there are $T',S'$, and $\beta \in \kappa$ such that 
\begin{enumerate}
\item \label{Hall}$\alpha = \omega_1 \cap \Hall(X, T', S', \beta \}; \Lsc_\kappa[X, T', S'])$ and 
\item \label{rest} $T \rest \alpha= T'\rest \alpha$ and $S \rest \alpha = S'\rest \alpha$,
\end{enumerate}
then $g(\alpha) = \otp(\Hall(\{ X, T', S', \beta \}; \Lsc_\kappa[X, T', S']) \cap \kappa)$. 
If there are no $T',S', \beta$ which satisfy the conditions above, 
let $g(\alpha)= f(\alpha)$.

We show that $g$ is well-defines as a function. 
In order to see this assume \ref{Hall} and \ref{rest} of 
the definition of $g$ hold for $(T^i, S^i, \beta^i)$, $i \in 2$.
Let $\gamma^i = \otp (\Hall(\{X, T^i, S^i, \beta \}); \Lsc_\kappa[X, T^i, S^i])$. Then by Lemma \ref{magic}, $\alpha$ is countable in both $\Lsc_{\gamma^i+2}[X \cap \alpha, T^i \rest \alpha, S^i \rest \alpha]$ 
for $i \in 2.$
But $T^0 \rest \alpha = T^1 \rest \alpha$
and $S^0 \rest \alpha = S^1 \rest \alpha $. 
Therefore, $\gamma^0$ and $ \gamma^1$ cannot be different limit ordinals.

Now we are ready to define $T(\alpha)$. 
If $b \subset T \rest \alpha$ is a cofinal branch, we let $\bigcup b \in T(\alpha)$ if and only if $b \in \Lsc_{f(\alpha)}[X \cap \alpha, T\rest \alpha, S\rest \alpha]$.
If $b \subset S \rest \alpha $ is a cofinal branch, we let 
$\bigcup b \in S(\alpha)$ if and only if 
$b \in \Lsc_{g(\alpha)} [X \cap \alpha, T \rest \alpha , S \rest \alpha]$.

In order to see $T$ does not club embed into $S$, 
assume for a contradiction that $C \subset \omega_1$ is a club and 
$k : T \rest C \longrightarrow S \rest C$ is a club embedding.
Since $X$ is $\kappa$-maximal, there is $\beta \in \kappa$, $\alpha \in \omega_1$ such that: 
\begin{itemize}
\item $\alpha = \omega_1 \cap \Hall(\{X, T, S, \beta \}; \Lsc_\kappa[X,T,S])$ and 
\item $\alpha = \omega_1 \cap \Hall(\{X, T, S,k, \beta \}; \Lsc_\kappa[X,T,S,k]).$
\end{itemize}
Then $g(\alpha)= \otp (\kappa \cap  \Hall(\{X, T, S, \beta \}; \Lsc_\kappa[X,T,S]))$.
Let $b$ be a cofinal branch of $T\rest \alpha$ which is in $ \Lsc_{f(\alpha)}[X \cap \alpha, T \rest \alpha, S\rest \alpha]$
such that $\alpha$ is countable in $\Lsc_{g(\alpha)}[X \cap \alpha, T \rest \alpha, S \rest \alpha, b]$.
For instance we can choose $b$ to be a cofinal branch of 
$\{t \in T \rest \alpha : \supp(t)$ is finite.$\}$ in such a way 
that $b$ codes a cofinal sequence of $\alpha$ with order type $\omega$.

Let $d = \bigcup k[b \rest C]$. Since $k$ is  a club embedding, $d \in S(\alpha)$ and consequently it is in 
$ \Lsc_{g(\alpha)}[X \cap \alpha, T \rest \alpha, S \rest \alpha, k \rest (T \rest (C \cap \alpha))].$
Let $\gamma = \kappa \cap \otp(\Hall(\{X, T , S, k, \beta \}; \Lsc_{\kappa}[X, T, S, k]))$. Observe $\gamma \geq g (\alpha)$.
Then $b =k^{-1}[d] \in \Lsc_\gamma [X \cap \alpha, T \rest \alpha, S \rest \alpha, k \rest (T \rest (C \cap \alpha))]$.
This is a contradiction because the last model 
thinks $\alpha$ is uncountable.

It suffices to show $S$ is Kurepa since it is a subtree of $T$.
Assume for a contradiction that $\Seq{b_\xi : \xi \in \omega_1}$
is the $\Lsc[X,S]$-minimal  enumeration of all cofinal branches of $S$.
Define $\Seq{M_\alpha: \alpha \in \omega_1}$ inductively as follows:
\begin{itemize}
\item $M_0 = \Hall(\{X,S\};\Lsc_{\omega_2}[X,S])$, 
\item $M_{\alpha+1} = \Hall(\{ M_\alpha \}; \Lsc_{\omega_2}[X,S])$, and 
\item $M_\alpha = \bigcup_{\xi \in \alpha}M_\xi$ if $\alpha \in \omega_1$ is limit.
\end{itemize} 
Let $\delta_\xi = M_\xi \cap \omega_1$, $\gamma_\xi = \otp(M_\xi \cap \omega_2)$, 
and 
$E = \{ \delta_\xi : \xi \in \omega_1 \}$.
Obviously, $E \subset \omega_1$ is a club and $M_\xi$ is 
isomorphic to $\Lsc_{\gamma_\xi}[X \cap \delta_\xi , T \rest \delta_\xi]$. 
Moreover $E \cap \delta_\xi \in \Lsc_{\gamma_\xi + 3} [X \cap \delta_\xi , T \rest \delta_\xi]$, as in Proposition \ref{Tree}. 

In order to reach a contradiction, 
we define inductively an increasing sequence $\Seq{t_\alpha : \alpha \in E}$ in $S$ such that $\bigcup_{\alpha \in E}t_\alpha$ is different from all branches in  $\Seq{b_\xi : \xi \in \omega_1 }$. 
Let $t_{\delta_0}$ be the $\Lsc[X,T]$-minimal element in $S(\delta_0 +1)$ which is not in $\bigcup_{\xi \in \delta_0} b_\xi$.
If $t_{\delta_\xi}$ is given let $t_{\delta_{\xi +1 }}$ be the 
$\Lsc[X,T]$-minimal element in $S(\delta_{\xi +1} +1)$ which is above 
$t_{\delta_\xi}$ and which is not in $\bigcup_{\eta \in \delta_{\xi +1}} b_\eta.$
If $\xi \in \omega_1$ is limit, let $t_{\delta_\xi} =  \bigcup_{\eta \in C \cap \delta_\xi}t_\eta$.
We claim that $t_{\delta_\xi}$ is defined for every $\xi \in \omega_1$. 
For the successor step this is obvious. 
For the limit $\xi \in \omega_1$, similar argument to the one in Proposition 
\ref{Tree} shows that $\gamma_\xi + \omega \leq g(\alpha)$, 
which implies that $t_{\delta_\xi}$ is defined.
\end{proof}

The following proposition completes the proof of Theorem \ref{main} and Theorem \ref{club_iso_theorem}.
\begin{prop}\label{max_in_semi}
Assume $\alpha_0 \in \omega_1$, $X \subset \omega_1$,
$\omega_1^{\Lsc[X\cap \alpha_0]} = \omega_1$, $\omega_2^{\Lsc[X]} = \omega_2$
$\mathbb{P}$ is a semiproper forcing notion in $\Lsc[X]$
such that $\textsc{V}$ is a forcing extension of $\Lsc[X]$ by $\mathbb{P}$.
Then $X$ is $\kappa$-maximal for all large enough regular cardinals $\kappa$. 
\end{prop}
\begin{proof}
Let $\kappa$ be a regular large enough cardinal in $\textsc{V}$ such that $2^\mathbb{P}$ is in $ \Lsc_\kappa [X]$ and 
whenever  $M \prec L_\kappa[X]$ is countable,  $\mathbb{P} \in M$
and $p \in M \cap \mathbb{P}$, then there is $q\leq p$ in $\mathbb{P}$ which is $(M, \mathbb{P})$-generic.

Assume for a contradiction that there is $p \in \mathbb{P}$ which forces that $X$ is not $\kappa$-maximal.
Let $\dot{Y}$ be a $\mathbb{P}$-name in $\Lsc_\kappa[X]$ for a subset of $\omega_1$ such that for all $\beta \in \kappa$, 
$$p \Vdash h_\kappa(X, \beta) < h_\kappa(\dot{Y}, \beta). $$
In particular, for all $\beta \in \kappa$,  $p$ forces that 
\begin{equation}
\omega_1 \cap \Hall(\{X, \beta \}, \Lsc_\kappa[X])
< \omega_1 \cap \Hall(\{X, \dot{Y}, \beta \}, \Lsc_\kappa[X, \dot{Y}]).
\end{equation}
Let $\beta$ be the order type of the predecessors of $(\mathbb{P}, p, \dot{Y}, X)$ in the well ordering of $\Lsc[X]$ and let $N = \Hall(\{ X, \beta \}; \Lsc_\kappa[X])$. Then $(\mathbb{P}, p, \dot{Y}, X)\in N.$ 
Let $\dot{G}$ be the canonical $\mathbb{P}$-name for the generic filter of $\mathbb{P}$.
Let $q \leq p$ such that $$q \Vdash N [\dot{G}] \cap \omega_1 = N \cap \omega_1.$$

On the other hand, since $q$ extends $p$ it forces that 
$$N \cap \omega_1 = \omega_1 \cap \Hall(\{X, \beta \}; \Lsc_\kappa[X]) < $$
$$< \omega_1 \cap \Hall(\{ X, \dot{Y}, \beta \}; \Lsc_\kappa[X,\dot{Y}]) \leq N[G] \cap \omega_1.$$
But this is  a contradiction since the first and the last ordinal 
are forced by $q$ to be equal.
\end{proof}

This proposition leads to the following natural questions.

\begin{question}
Suppose that the constructible universe $L$ has no inaccessible cardinals. Is there $X\subset \omega_1$ such that $\omega_1^{\Lsc[X\cap \alpha_0]} = \omega_1$, $\omega_2^{\Lsc[X]} = \omega_2$ for some $\alpha_0\in \omega_1$ and such that $X$ is $\kappa$-maximal for all large enough regular cardinals $\kappa$? What if we assume further that $V$ is a forcing extension of $L[X]$ that preserves stationary subsets of $\omega_1$?
\end{question}

Let  $X\subset \omega_1$ such that $\omega_1^{\Lsc[X\cap \alpha_0]} = \omega_1$, $\omega_2^{\Lsc[X]} = \omega_2$ for some $\alpha_0\in \omega_1$. Let $\rho \in L[X]$ be a subadditive function (see \cite{walks}) mapping $[\omega_2]^2$ to $\omega_1.$  
Let $f_\alpha (\nu) = \otp (\{ \xi \in \alpha: \rho(\xi , \alpha) < \nu \})$,
for $\alpha \in \omega_2$ and $\nu \in \omega_1$. It is easy to see that 
$\Seq{f_\alpha: \omega_1 \longrightarrow \omega_1: \alpha \in \omega_2}$ is strictly increasing modulo the ideal of countable sets.
It is also easy to see that the function $f$ of the proof of 
Proposition \ref{Similar} dominates all the functions of any such sequence.
This inspires the following question, 
which can also be viewed as an idea for exploring the relation
between Aronszajn free Kurepa trees and the consistency of the existence of inaccessible cardinals.  

\begin{question}
Suppose that the constructible universe $L$ has no inaccessible cardinals. Assume for every sequence of 
functions 
$\Seq{f_\alpha : \omega_1 \longrightarrow \omega_1 : \alpha \in \omega_2}$ there is $g: \omega_1 \longrightarrow \omega_1$
which dominates all 
$f_\alpha$'s on a club.
Does there exist 
$X \subset \omega_1$ and $\alpha_0 \in \omega_1$
such that $\omega_1= \omega_1^{\Lsc[X \cap \alpha_0]}$, $\omega_2 = \omega_2^{\Lsc[X]}$ and $X$ is $\kappa$-maximal for all large enough regular $\kappa$?
\end{question}

\vspace{3mm}
Acknowledgement. The research on this paper was partially supported by grants from NSERC(455916) and SFRS(7750027-SMART).

\def\Dbar{\leavevmode\lower.6ex\hbox to 0pt{\hskip-.23ex \accent"16\hss}D}

\end{document}